\documentclass[11pt]{amsart}
\usepackage{amssymb}

\newcommand{\Lie}{\mathcal}
\renewcommand{\sl}{\operatorname{sl}}
\newcommand{\GL}{\operatorname{GL}}
\newcommand{\Ad}{\operatorname{Ad}}
\newcommand{\Rad}{\operatorname{Rad}}
\newcommand{\nil}{\operatorname{nil}}
\newcommand{\iso}{\cong}
\newcommand{\prodsemi}{\rtimes}
\newcommand{\integer}{\mathord{\mathbb Z}}
\newcommand{\normal}{\triangleleft}
\newcommand{\real}{\mathord{\mathbb R}}
\newcommand{\func}{\operatorname{Func}}

\newcommand{\cover}{\widetilde}

\newcommand{\ma}{\mathcal{M}}
\newcommand{\ka}{\mathcal{K}}
\newcommand{\sa}{\mathcal{S}}
\newcommand{\pa}{\mathcal{P}}
\newcommand{\fa}{\mathcal{F}}
\newcommand{\kq}{{\overline{\ka}}}
\newcommand{\aq}{{\overline{\alpha}}}

\newcommand{\mub}{\overline{\mu}}
\newcommand{\psib}{{\overline{\psi}}}

\newcommand\set{} 
\def\set\{#1\mid #2\}{\left\{ #1
 \left| \vphantom {\left\{ #1 \mid #2 \right\} } \right.
 #2 \right\} }

\renewcommand{\see}[1]{\textnormal{(}see~\protect\ref{#1}\textnormal{)}}
\newcommand{\seemore}[2]{\textnormal{(}see
#1~\protect\ref{#2}\textnormal{)}}
\newcommand{\cf}[1]{(cf.~\protect\ref{#1})}
\newcommand{\pref}[1]{\textnormal{(\protect\ref{#1})}}

\newcommand \step[2]{\medskip \noindent
 \emph{Step \textnormal{#1.} #2} \ignorespaces}

\renewcommand{\thmhead}[3]{%
 \thmname{#1}\thmnumber{ #2}\thmnote{ {\the\theoremnotefont#3}}}

\renewcommand{\swappedhead}[3]{%
 \textnormal{(\thmnumber{#2})}\thmname{ #1}\thmnote{
 {\the\theoremnotefont#3}}}
\swapnumbers

\newtheorem{thm}{Theorem}[section]
\newtheorem{lem}[thm]{Lemma}
\newtheorem{cor}[thm]{Corollary}
\newtheorem{prop}[thm]{Proposition}

\theoremstyle{definition}

\newtheorem{notation}[thm]{Notation}
\newtheorem{standing}[thm]{Standing assumptions}
\newtheorem{defn}[thm]{Definition}
\newtheorem{assump}[thm]{Assumption}

\theoremstyle{remark}

\newtheorem{rem}[thm]{Remark}        \renewcommand{\therem}{}
\newtheorem{ack}[thm]{Acknowledgment}  \renewcommand{\theack}{}

\newbox\sectionS \setbox\sectionS\hbox{\S}
\renewcommand{\thesubsection}{\copy\sectionS
 \thesection\Alph{subsection}}

\newcommand{\refjour}{} 
\newcommand{\refbook}{}
\newcommand{\refappear}{}

\long\def\refjour[#1]#2 #3: #4. #5.. #6 (#7) #8 \par
 {\bibitem{#2} #3, #4, \emph{#5} {\bf #6} (#7), #8.\par}
\long\def\refbook[#1]#2 #3: #4. #5, #6, #7 \par
 {\bibitem{#2} #3, ``#4,'' #5, #6, #7.\par}
\long\def\refappear[#1]#2 #3: #4. #5 \par
 {\bibitem{#2} #3, #4, \emph{#5} (to appear).\par}

\begin{document}

\title[Cocycle superrigidity for non-semisimple groups]
{Cocycle superrigidity for ergodic actions \\
 of non-semisimple Lie groups}

\author{Dave Witte}

\address{Department of Mathematics, Oklahoma State University,
Stillwater, OK 74078}

\email{dwitte@math.okstate.edu}

\date{\today}

\thanks{Submitted in June 1996 to the \emph{Proceedings of the
Colloquium on Lie Groups and Ergodic Theory
(Bombay, India, 4--12 January 1996)}.}

\begin{abstract}
 Suppose $L$ is a semisimple Levi subgroup of a connected Lie group~$G$,
$X$ is a Borel $G$-space with finite invariant measure, and $\alpha
\colon X \times G \to \GL_n(\real)$ is a Borel cocycle. Assume $L$ has
finite center, and that the real rank of every simple factor of~$L$ is
at least two. We show that if $L$ is ergodic on~$X$, and the restriction
of~$\alpha$ to~$X \times L$ is cohomologous to a homomorphism (modulo a
compact group), then, after passing to a finite cover of~$X$, the
cocycle $\alpha$ itself is cohomologous to a homomorphism (modulo a
compact group).
 \end{abstract}

\maketitle

\section{Introduction} \label{introduction}

\begin{defn}[{\cite[Defns.~4.2.1 and 4.2.2, p.~65]{ZimmerBook}}]
\label{cocycle-def}
 Suppose $G$ and~$H$ are Lie groups, and $X$ is a Borel $G$-space. A
Borel function $\alpha \colon X \times G \to H$ is a \emph{Borel
cocycle} if, for all $g,h \in G$, we have
 $$ (x,gh)^\alpha = (x,g)^\alpha (xg,h)^\alpha
 \qquad \textnormal{for a.e.~$x \in X$} . $$
 Two cocycles $\alpha $~and~$\beta$ are \emph{cohomologous} if there
is a Borel function $\phi \colon X \to H$, such that, for all $g \in
G$, we have
 $$ (x,g)^\alpha = (x^\phi )^{-1} (x,g)^\beta (xg)^\phi
 \qquad \textnormal{for a.e.~$x \in X$} . $$
 \end{defn}

Any continuous group homomorphism $\sigma \colon G \to H$ gives rise
to a cocycle, defined by $( \cdot ,g)^\alpha = g^\sigma $. For actions
of semisimple groups, R.~J.~Zimmer's Cocycle Superrigidity Theorem often
shows that (up to cohomology, and modulo a compact group) these are the
only examples.

\begin{defn}[{(cf.~\ref{Zar-hull-def})}]
 Suppose $G$ is a Lie group, $X$ is an ergodic Borel $G$-space with
finite invariant measure, $H$~is a subgroup of $\GL_n(\real)$, and
$\alpha \colon X \times G \to H$ is a Borel cocycle. We say $\alpha$ is
\emph{Zariski dense} if $H$~is contained in the Zariski closure of the
range of every cocycle $\beta \colon X \times G\to H$ that is
cohomologous to~$\alpha $.
 \end{defn}

\begin{thm}[{(Zimmer, cf.~\cite[Thms.~5.2.5, 7.1.4, 9.1.1]{ZimmerBook})}]
\label{Zimmer-super}
 Suppose $G$ is a connected, semisimple Lie group, $X$ is an ergodic
Borel $G$-space with finite invariant measure, $H$~is a Zariski closed
subgroup of $\GL_n(\real)$, and $\alpha \colon X \times G \to H$ is a
Zariski-dense Borel cocycle. Assume $G$ has finite center, and that the
real rank of every simple factor of~$G$ is at least two. If $H$ is
reductive, then, after replacing~$G$ and~$X$ by finite covers, there are:
 \begin{itemize}
 \item a homomorphism $\sigma \colon G \to H$; 
 \item a compact, normal subgroup~$K$ of $H$ that centralizes~$G^\sigma
$; and
 \item a cocycle~$\beta$ that is cohomologous to~$\alpha $;
 \end{itemize}
 such that, for every $g \in G$, we have $(x,g)^\beta \in g^\sigma K$
for a.e.~$x\in X$.
 \end{thm}

This paper extends Zimmer's result to groups that are not semisimple.
Our main theorem reduces the general case to the semisimple case.

\begin{thm} \label{mainthm}
 Assume
 \begin{itemize}
 \item $G$ is a connected Lie group;
 \item $X$ is a Borel $G$-space with finite invariant measure;
 \item $H$ is a connected Lie subgroup of~$\GL_n(\real)$ that is of
finite index in its Zariski closure, and has no nontrivial compact,
normal subgroups;
 \item $\alpha \colon X \times G \to H$ is a Zariski-dense Borel cocycle;
 \item $L$ is the product of the noncompact, simple factors in a
semisimple Levi subgroup of~$G$;
 \item $L$ is ergodic on~$X$; 
 \item $\sigma \colon L \to H$ is a continuous homomorphism;
 \item $K$ is a compact subgroup of~$H$ that centralizes~$L^\sigma $;
 \item $H = (\Rad H) L^\sigma K$;
 and
 \item for every $l \in L$, we have $(x,l)^\alpha \in l^\sigma K$ for
a.e.~$x \in X$.
 \end{itemize}
 Then
 \begin{itemize}
 \item $\sigma$ extends to a continuous homomorphism defined on all
of~$G$; and
 \item $\alpha$ is cohomologous to the cocycle~$\beta$, defined by
$( \cdot ,g)^\beta = g^\sigma$.
 \end{itemize}
 \end{thm}

By combining our theorem with Zimmer's, we obtain the following general
result.

\begin{defn}[{(cf. \cite[Defn.~9.2.2, p.~167]{ZimmerBook})}]
\label{Zar-hull-def}
 Suppose $L$ is a connected Lie group, $X$ is an ergodic Borel
$L$-space with finite invariant measure, and $\alpha \colon X \times L
\to \GL_n(\real)$ is a Borel cocycle. The \emph{Zariski hull} of~$\alpha
$ is a Zariski closed subgroup~$J$ of $\GL_n(\real)$, such that $\alpha$
is cohomologous to a cocycle $\beta \colon X \times L \to
\GL_n(\real)$, such that the range of~$\beta$ is contained in~$J$, but
$\alpha$ is not cohomologous to any cocycle whose range is contained in
a Zariski-closed, proper subgroup of~$J$. The Zariski hull~$J$ always
exists, and is unique up to conjugacy.
 \end{defn}

\begin{cor} \label{main-cor}
 Assume
 \begin{itemize}
 \item $G$ is a connected Lie group;
 \item $X$ is a Borel $G$-space with finite invariant measure;
 \item $H$ is a Lie subgroup of~$\GL_n(\real)$ that
is of finite index in its Zariski closure, and has no nontrivial compact,
normal subgroups;
 \item $\alpha \colon X \times G \to H$ is a Zariski-dense Borel cocycle;
 \item $L$ is the product of the noncompact, simple factors in a
semisimple Levi subgroup of~$G$;
 \item $L$ is ergodic on~$X$; 
 \item $L$ has finite center, and the real rank of every simple factor
of~$L$ is at least two; and
 \item the Zariski hull of the restriction of~$\alpha$ to $X \times L$
is reductive.
 \end{itemize}
 Then, after passing to a finite cover of~$X$, the cocycle~$\alpha$ is
cohomologous to a homomorphism.
 \end{cor}

\begin{rem}
 The assumption that $L$ is ergodic on~$X$ cannot be weakened to the
ergodicity of~$G$. To see this, suppose $L$ is not ergodic, and let $Y$
be the space of ergodic components of~$L$ on~$X$. The Mautner
phenomenon \cite[Thm.~1.1]{Moore-Mautner} implies that $L$ and $[L,G]$
have the same ergodic components, so the $G$-action on~$X$ factors
through to an action of $G/[L,G]$ on~$Y$. Because $G/[L,G]$ is amenable,
there may be cocycles of this action that are not related to the
algebraic structure of the acting group (cf.~\cite{ConnesFeldmanWeiss}).
Pulling back to~$G$, these are cocycles on~$X$ that have almost nothing
to do with~$G$.

On the other hand, it is not known whether the assumption that $H$ is
reductive can be omitted from Zimmer's Theorem; perhaps this hypothesis
is always satisfied. Zimmer \cite[cf.~Thm.~1.1]{Zimmer-reductive}
proved this to be the case for cocycles that satisfy an $L^1$~growth
condition. Then, by an argument very much in the spirit of the present
paper, he was able to derive a superrigidity theorem for $L^1$~cocycles
of actions of some non-semisimple groups (see
\cite[Thm.~4.1]{Zimmer-reductive}).
 \end{rem}

 Zimmer's cocycle superrigidity theorem~\pref{Zimmer-super} was inspired
by the superrigidity theorem for finite-dimensional representations of
lattices in semisimple Lie groups, proved by G.~A.~Margulis
\cite[Thm.~VII.5.9, p.~230]{MargBook}. The present work was suggested by
the author's \cite{Witte-super}, \cite[\S2, \S5]{Witte-super-S}
generalization of Margulis' theorem to non-semisimple groups.

After some preliminaries in~\S\ref{prelims}, we prove a restricted version
of Thm.~\ref{mainthm} in~\S\ref{main-proof}. The final section of the
paper removes the restrictions, and presents a proof of
Cor.~\ref{main-cor}.

\begin{ack}
 I would like to thank my colleagues at the Tata Institute of Fundamental
Research for their generous hospitality throughout the visit during
which much of this work was carried out. I would also like to thank Alex
Eskin and Robert J.~Zimmer for several helpful suggestions. This
research was partially supported by a grant from the National Science
Foundation.
 \end{ack}

\section{Preliminaries} \label{prelims}

\begin{standing}
The notation and hypotheses of Thm.~\ref{mainthm} are in effect
throughout this section.
 \end{standing}

\begin{notation}
 We usually write our maps as superscripts. Thus, if $x\in X$ and $\phi 
\colon X \to Y$, then $x^\phi$ denotes the image of~$x$ under~$\phi $.

If $g$~and~$h$ belong to~$H$ (or to any other group), then $g^h$
denotes the conjugate $h^{-1} gh$.

 For any Borel function~$\phi$ whose range lies in~$H$, we use $-\phi$
to denote the function defined by
 $( \cdot )^{-\phi } = \bigl( ( \cdot )^\phi \bigr)^{-1}$.
 \end{notation}

\begin{defn}
 Let $Q$ be a subset of~$G$, and let $\fa \colon X \times G \to H$ be
a Borel function. We use $\fa|_Q$ to denote the restriction of~$\fa$ to
$X \times Q$.

If $Q$ is countable, or is a Lie subgroup of~$G$, then there is a
natural choice of a measure class on~$Q$, and we use $(X \times Q)^\fa$
to denote the essential range of~$\fa|_Q$. (Recall that the
\emph{essential range} is the unique smallest closed set whose inverse
image is conull.)

If, for some $r \in G$, the function $\fa|_r$ is essentially constant,
then we often omit the reference to~$X$, and simply write~$r^\fa$ for
the single point in $(X \times r)^\fa$.
 \end{defn}

\begin{defn}
 Let $Q$ be a subgroup of~$G$, and let $\fa \colon X \times G \to H$ be
a Borel function.
 We say that $\fa|_Q$ is a \emph{homomorphism} if there is a
homomorphism $\sigma \colon Q \to S$, such that, for all $r \in Q$, we
have $( \cdot ,r)^\fa = r^\sigma$~a.e.
 \end{defn}

The following well-known result is a straightforward consequence of the
cocycle identity.

\begin{lem} \label{homo-const}
 Let $Q$ be a subgroup of~$G$, and let $\fa \colon X \times G \to H$ be
a Borel cocycle. Then $\fa|_Q$ is a homomorphism iff $( \cdot ,r)^\fa
= r^\fa$ is essentially constant, for each $r \in Q$. \qed
 \end{lem}

\begin{defn}
 Let us say that an element $g \in \GL_n(\real)$ is \emph{split} if every
eigenvalue of~$g$ is real and positive. (In other words, the real Jordan
decomposition of~$g$ \cite[Lem.~IX.7.1, p.~430]{Helgason} has no
elliptic part.)

Let us say that an element of~$L$ is \emph{split} if it belongs to a
one-parameter subgroup~$T$ of~$L$, such that $\Ad_L t$ is split, for all
$t \in T$. Note that if $l$~is a split element of~$L$, then $L^\sigma$ 
is a split element of~$H$.
 \end{defn}

\begin{lem} \label{must-norm}
 Suppose $l$ is a split element of~$L$, $k \in K$, and $T$ is
a connected subgroup of~$H$, such that $l^\sigma k$ normalizes~$T$. Then
$k$ normalizes~$T$.
 \end{lem}

\begin{proof}
 The normalizer of~$T$ is Zariski closed (cf.~pf{.} of \cite[Thm.~3.2.5,
p.~42]{ZimmerBook}), so it contains the elliptic part of each of its
elements (cf.~proof of \cite[Thm.~15.3, p.~99]{Humphreys-alg}).
Therefore, the desired conclusion follows from the observation that $k$
is the elliptic part of~$l^\sigma k$.
 \end{proof}

\begin{thm}[{(Borel Density Thm. \cite[Cor.~2.6]{Dani-BDT})}]
\label{BDT}
 Suppose $H$ acts regularly on a variety~$V$, $\mu$ is a
probability measure on~$V$, and $h$ is a split element of~$H$. If $\mu$
is $h$-invariant, then $h$ fixes the support of~$\mu$ pointwise. \qed
 \end{thm}

\begin{cor} \label{BDT-cor}
 Suppose $H$ acts regularly on a variety~$V$.
 Let $\psi \colon X \to V$ be a Borel function, and let $g \in G$ and
$h \in H$, and let $K$ be a compact subgroup of~$H$ that
centralizes~$h$. Assume $x^{g\psi} \in (x^{\psi h})^K$ for a.e.~$x \in
X$, and $h$~is split.
 Then $h$ fixes the essential range of~$\psi$ pointwise.
 \end{cor}

\begin{proof}
 Let $V/K$ be the space of $K$-orbits on~$V$. The induced map $\psib
\colon X \to V/K$ satisfies $x^{g\psib} = x^{\psib h}$. Therefore, the
$G$-invariant probability measure on~$X$ pushes via~$\psib$ to an
$h$-invariant probability measure~$\mub$ on $V/K$.
By using~$\mub$ to integrate together the $K$-invariant probability
measure on each $K$-orbit, we may lift~$\mub$ to an
$h$-invariant probability measure~$\mu$ on~$V$.
 Note that the support of~$\mu$ is $(X^\psi)^K$, where $X^\psi$ is the
essential range of~$\psi$. On the other hand, the Borel Density Theorem
implies that $h$ fixes the support of~$\mu$ pointwise.
 \end{proof}

\begin{thm}[{(Moore Ergodicity Theorem, cf.~\cite[Thm.~1.1]{Moore-Mautner})}]
\label{MooreErgodicity}
 If $T$ is a connected subgroup of~$L$, such that, for every simple
factor~$L_i$ of~$L$, the projection of~$T$ into $L_i/Z(L_i)$ has
noncompact closure, then $T$ is ergodic on~$X$. \qed
 \end{thm}

\begin{prop}[{\cite[Thm.~XV.3.1, pp.~180--181, and see
p.~186]{Hochschild-Lie}}] \label{cpct-conj}
 If $G$ is a Lie group that has only finitely many connected
components, then $G$ has a maximal compact subgroup~$K$, and every
compact subgroup of~$G$ is contained in a conjugate of~$K$. \qed
 \end{prop}

\begin{prop}[{\cite[Thm.~XIII.1.3, p.~144]{Hochschild-Lie}}]
\label{cpct-abel} 
 If a Lie group~$G$ is compact, connected, and solvable, then $G$ is
abelian. \qed
 \end{prop}

\begin{lem} \label{cover-coho}
 Let $\cover{G}$ be a covering group of~$G$, and let $\cover{\alpha }
\colon X \times \cover{G} \to H$ be the cocycle naturally induced
by~$\alpha $. If $\cover{\alpha}$ is cohomologous to a homomorphism,
then $\alpha$ is cohomologous to a homomorphism.
 \end{lem}

\begin{proof}
 By assumption, there is a Borel function $\phi \colon X \to H$, such
that, for all $g \in \cover{G}$, the expression
 $x^{-\phi } (x,g)^{\cover{\alpha }} (xg)^\phi$
 is essentially independent of~$x$. Then the same is true with $\alpha$
in place of~$\cover{\alpha }$, for all $g \in G$, as desired.
 \end{proof}

\section{Proof of Theorem~\ref{mainthm} (The Main Case)}
\label{main-proof}

This entire section is devoted to a proof of Thm.~\ref{mainthm}, so the
notation and hypotheses of Thm~\ref{mainthm} are in effect throughout.

\begin{assump} \label{special-case}
 Throughout this section, we assume $\Rad G$ is nilpotent, and that $G$
has no nontrivial compact semisimple quotients. See \ref{general-case}
for an explanation of how to obtain the full theorem from this special
case.
 \end{assump}

\begin{notation}
 Let $R = \Rad G$, so $G = RL$. By passing to a covering group of~$G$, we
may assume $R$ is simply connected \see{cover-coho}.

By assumption (and perhaps replacing~$K$ with a larger compact group), we
may write $H = S \prodsemi (MK)$, where
 \begin{itemize}
 \item $S$ is a connected, split, solvable subgroup;
 \item $M = L^\sigma$ is connected and semisimple, with no compact
factors;
 \item $K$ is a compact subgroup that centralizes~$M$;
 \item $(X \times L)^\alpha \subset MK$; and
 \item $M \cap K = Z(M)$.
 \end{itemize}
 \end{notation}

Now $\alpha$ induces a cocycle $\aq \colon X \times G \to H/(SK)
\iso M/(M \cap K)$. Note that $\aq|_L = \overline{\sigma }$ is a
homomorphism.

\begin{prop}
 $(X \times R)^\aq = e$.
 \end{prop}

\begin{proof}
 Let $P$ be a minimal parabolic subgroup of~$G$. From
 \cite[Step~1 of pf.~of Thm.~5.2.5, p.~103]{ZimmerBook},
 we know there is an (almost) Zariski closed, proper subgroup~$L$ of~$M$
and a Borel function $\phi \colon X \times P \backslash G \to
L\backslash M$ such that, for all $g \in G$, we have
 $$ (xg, cg)^\phi = (x,c)^\phi (x,g)^\aq
 \qquad \textnormal{for a.e.~$(x,c) \in X \times P\backslash G$} .$$
 In particular, for $l \in L$, we have
 $(xl, cl)^\phi = (x,c)^\phi l^\aq$.
 Then, by Fubini's Theorem, we see that, for a.e.~$c \in P\backslash
G$, $l \in P^c \cap L$, and $x \in X$, we have
 $(xl, c)^\phi = (x,c)^\phi l^\aq$.
 So \ref{BDT-cor} implies that $l^\aq$ fixes $(X \times c)^\phi$
pointwise (assuming that $l$~is split), which implies $(xl,c)^\phi =
(x,c)^\phi $.
 So, from the ergodicity of~$P^c \cap L$ \see{MooreErgodicity}, we
conclude that $( \cdot ,c)^\phi = c^\phi$ is essentially constant.
Therefore, for $c \in P\backslash G$ and $r \in \Rad G$, because $cr =
c$, we have
 $$ (x,c)^\phi = (xr, cr)^\phi = (x,c)^\phi (x,r)^\aq .$$
 Therefore, $(x,r)^\aq$ fixes $(X \times P\backslash G)^\phi $
pointwise, so $(x,r)^\aq$ is trivial (cf.~\cite[pf.~of Lem.~5.2.8,
p.~102]{ZimmerBook}).
 \end{proof}

\begin{defn}
 For $r \in R$ and $l \in L$, we have
 $(x,rl)^\aq = (x,r)^\aq (xr,l)^\aq = l^\aq$, so we see that $\aq$ is a
homomorphism.
 By replacing $G$ with a finite cover, we may assume that $\aq$ lifts to
a homomorphism $\ma\colon G \to M$ \see{cover-coho}. Because $H = S
\prodsemi (MK)$, there are well-defined Borel functions $\sa \colon X
\times G \to S$ and $\ka \colon X \times G \to K$, such that
 $$(x,g)^\alpha = (x,g)^\sa g^\ma (x,g)^\ka
 \qquad \textnormal{for a.e.~$x \in X$.}$$
 \end{defn}

Note that, for $u \in L$ and $r \in R$, we have
 \begin{itemize}
 \item $(x,u)^\alpha = u^\ma (x,u)^\ka$ for a.e.~$x \in X$; and
 \item $(x,r)^\alpha = (x,u)^\sa (x,u)^\ka$ for a.e.~$x \in X$.
 \end{itemize}
Note also that $G^\ma$ centralizes $(X \times G)^\ka$,
because $M$ centralizes~$K$.

Because $\sa$ may not be a cocycle, Lem.~\ref{homo-const} may not
apply to~$\sa$. However, the following lemma is a suitable replacement.

\begin{defn}
 Let $Q$ be a subgroup of~$R$. A function $\sigma \colon Q \to S$ is a
\emph{crossed homomorphism} if there is a homomorphism $\kappa \colon R
\to K'$, such that $(rs)^\sigma = r^\sigma s^{\sigma r^\kappa}$ for all
$r,s \in Q$, where
 $K' = N_K( \langle Q^\sigma \rangle )/C_K( \langle Q^\sigma \rangle
)$.

 We say that $\sa|_Q$ is a \emph{crossed homomorphism} if there is a
crossed homomorphism $\sigma \colon Q \to S$ such that, for all $r \in
Q$, we have $( \cdot ,r)^\sa = r^\sigma$~a.e.
 \end{defn}

\begin{lem} \label{when-crossed}
 Let $Q$ be a subgroup of~$R$.
 \begin{enumerate}
 \item \label{crossed-const}
 $\sa|_Q$ is a crossed homomorphism iff $( \cdot ,r)^\sa$ is essentially
constant, for each $r \in Q$.
 \item \label{homo-cent}
 $\sa|_Q$ is a homomorphism iff $\sa|_Q$ is a crossed homomorphism and
$(X \times Q)^\ka$ centralizes~$Q^\sa$.
 \end{enumerate}
 \end{lem}

\begin{proof} \pref{crossed-const} We need only prove the nontrivial
direction, so assume $( \cdot ,r)^\sa = r^\sa$ is essentially
constant, for each $r \in Q$.
 For convenience, let $N = N_K(\langle Q^\sa \rangle)$ and $C =
C_K(\langle Q^\sa \rangle)$.

 For $r,s \in Q$, we have 
 $(x,rs)^\alpha =(x,r)^\alpha (xr,s)^\alpha $,
 so $(rs)^\sa = r^\sa s^{\sa (x,r)^{-\ka}}$.
 This implies that $s^{\sa (x,r)^{-\ka}} \in \langle Q^\sa \rangle$,
so we see that $(X \times Q)^\ka \subset N$. This also implies that
$(\cdot,r)^\ka$ is essentially constant, modulo~$C$. Therefore, the
induced cocycle $\kq \colon X \times N/C$ is a homomorphism, as desired.

\pref{homo-cent} Because $(rs)^\sa = r^\sa s^{\sa (x,r)^{-\ka}}$, we
see that $(rs)^\sa = r^\sa s^\sa$ iff $(x,r)^\ka$ centralizes~$s^\sa$.
 \end{proof}

\begin{cor}[(of proof)] \label{crossed-K-homo}
 Suppose $Q$ is a subgroup of~$R$. Let
 $$N = N_K(\langle Q^\sa \rangle)
 \text{ and }
 C = C_K(\langle Q^\sa \rangle) .$$
 If $\sa|_Q$ is a
crossed homomorphism, then the cocycle $\kq \colon X \times Q\to N/C$,
induced by~$\ka$, is a homomorphism. \qed
 \end{cor}

\begin{notation}
 Let $U$ be a maximal connected unipotent subgroup of~$L$, let $U^-$ be
a maximal connected unipotent subgroup that is opposite to~$U$, and let
$A$ be a maximal split torus of~$L$ that normalizes both $U$ and~$U^-$.
Thus, $N_L(U) \cap N_L(U^-)$ is reductive, and contains~$A$ in its
center.
 \end{notation}

\begin{lem} \label{u-equivariant}
 Suppose $r \in R$, and $u$ is a split element of~$L$. Let $w = [u,r]
\in R$, and assume $( \cdot ,w)^\sa = w^\sa$ is essentially constant,
and that $(X \times u)^\ka$ and $(X \times w)^\ka$ centralize~$w^\sa$.
Then $(xu,r)^\sa = (x,r)^{\sa (x,u)^\ka}$ for a.e.~$x \in X$.
 \end{lem}

\begin{proof}
 For a.e.~$x \in X$, because $r = u^{-1}ruw$, we have
 \begin{eqnarray*}
 (xu,r)^\alpha
 &=& (x,u)^{-\alpha } (x,r)^\alpha (xr,u)^\alpha (xru, w)^\alpha \\
 &=& (x,u)^{-\ka} (x,r)^{\alpha u^\ma} w^\sa (xr,u)^\ka (xru, w)^\ka
\\
 &\in& \bigl( (x,r)^{\alpha u^\ma} w^\sa \bigr)^{(x,u)^\ka} C_K(w^\sa)
 \end{eqnarray*}
 So the Borel Density Theorem \see{BDT-cor} implies that $(x,r)^{\alpha
u^\ma} w^\sa = (x,r)^\alpha $. Thus, we have
 $$(xu,r)^\alpha \in (x,r)^{\alpha(x,u)^\ka}
C_K(w^\sa) ,$$
 which implies $(xu,r)^\sa = (x,r)^{\sa {(x,u)^\ka}}$.
 \end{proof}

\begin{cor} \label{Ku-cent}
 Suppose $r \in R$, $u$ is a split element of~$L$, and $W$ is a
subgroup of~$R$. Assume $\sa|_W$ is a homomorphism, $(X \times u)^\ka$
centralizes $W^\sa$, and that $w = [u,r] \in W$. If $( \cdot ,r)^\sa =
r^\sa$ is essentially constant, then $( X \times u)^\ka$
centralizes~$r^\sa$. \qed
 \end{cor}

\begin{cor} \label{Kunip-cent-W}
 Suppose $u$ is a unipotent element of~$L$, and $W$ is a
subgroup of~$R$ that is normalized by~$u$. If $\sa|_W$ is a
homomorphism, then $( X \times u)^\ka$ centralizes~$W^\sa$.
 \end{cor}

\begin{proof}
 Because $u$ is unipotent, we may assume by induction on $\dim W$ that
$(X \times u)^\ka$ centralizes $([u,W] [W,W])^\sa$. Then \ref{Ku-cent}
implies that $( X \times u)^\ka$ centralizes~$W^\sa$.
 \end{proof}

\begin{defn}
 Given a Borel function $\phi \colon X \to K$, let $\alpha ^\phi
\colon X \times G\to H$ be the cocycle cohomologous to~$\alpha $
defined by
 $$ (x,g)^{\alpha ^\phi } = x^{-\phi } (x,g)^\alpha (xg)^\phi .$$
 Also define a Borel function $\sa^\phi \colon X \times G \to \sa$ by
 $(x,g)^{\sa^\phi } = (x,g)^{\sa x^\phi }$.
 Note that
 $(x,g)^{\alpha ^\phi } \in (x,g)^{\sa^\phi } g^\ma K$,
 for all $(x,g) \in X \times G$.
 \end{defn}

\begin{cor} \label{wolg-const}
 Suppose $Q$~and~$W$ are Lie subgroups of~$R$, and $u$~is a split
element of~$L$, such that $[u,Q] \subset W$. Assume $\sa|_W$ is a
homomorphism, and that $(X \times u)^\ka$ centralizes~$W^\sa$. 
 Let $K_u$ be the closure of\/ $\langle ( X \times u)^\ka \rangle$. If
$u$ is ergodic on~$X$, then there is a Borel function $\phi \colon X
\to K_u$, such that $\sa^\phi |_Q$ is a crossed homomorphism.
 \end{cor}

\begin{proof}
 Let $\func(Q,S)$ be the the space of
Borel functions from~$Q$ to~$S$, where two functions are identified if
they agree almost anywhere. (The topology of convergence in measure
defines a countably generated Borel structure on this space
\cite[pp.~49--50]{ZimmerBook}.)
 The function $\sa|_Q \colon X \times Q \to S$ determines a Borel
function $\fa \colon X \to \func(Q,S)$. From \ref{u-equivariant}, we see
that, for a.e.~$x \in X$, we have
 $(xu)^\fa = x^{\fa (x,u)^\ka} \in x^{\fa K_u}$,
 where $K_u$ acts on $\func(Q,S)$ via conjugation on the range space~$S$.
Thus, the ergodicity of~$u$ implies that there is a Borel function
$\phi \colon X \to K_u$, and some $\sigma \in \func(Q,S)$, such that
 $x^{\fa x^\phi } = \sigma$ for a.e.~$x \in X$.
 That is, for a.e.~$x \in X$ and a.e.~$r \in Q$, we have
 $(x,r)^{\sa x^\phi } = r^\sigma $.
 In other words, $( \cdot ,r)^{\sa^\phi} = r^\sigma$ is essentially
constant, for a.e.~$r \in Q$. Then, because $Q$ has no proper subgroups
of full measure, we conclude from the cocycle identity (applied
to~$\alpha ^\phi$) that we have $( \cdot ,r)^{\sa^\phi} = r^\sigma$, for
\emph{all} $r \in Q$.
 \end{proof}

The following is the special case where $W$ is trivial.

\begin{cor} \label{wolg-crossed}
 Let $u$ be a split element of~$L$ that is ergodic on~$X$, and let
$K_u$ be the closure of $\langle ( X \times u)^\ka \rangle$. Then there
is a Borel function $\phi \colon X \to K_u$, such that $\sa^\phi
|_{C_R(u)}$ is a crossed homomorphism. \qed
 \end{cor}

Most of the work in this section is devoted to showing that we may
assume $\sa|_R$ is a homomorphism. The following proposition represents
our first real progress toward this goal. Most of the rest is achieved
by an inductive argument based on the unipotence of~$U$ and the
solvability of~$R$.

\begin{prop} \label{CRL-homo}
 We may assume $\sa|_{C_R(L)}$ is a homomorphism.
 \end{prop}

\begin{proof}
 For convenience, let $Q = C_R(L)$, let $V$ be the Zariski closure of
$\langle Q^\sa \rangle$, $N = N_K(V)$ and $C = C_K(V)$. From
\ref{wolg-crossed}, we see that, by replacing~$\alpha$ with a
cohomologous cocycle~$\alpha ^\phi$, we may assume $\sa|_{C_R(L)}$ is a
crossed homomorphism. Then \ref{crossed-K-homo} implies that the induced
cocycle $\kq \colon X \times Q\to N/C$ is a homomorphism. Therefore, the
induced cocycle $\aq \colon X \times Q \to VN/C$ is a homomorphism.
Then, because $Q$ is nilpotent \see{special-case}, the Zariski
closure of~$Q^\aq$ in $VN/C$ is of the form $W \times T/C$, where $W
\subset V$ is split and $T$ is a compact torus (cf.~\cite[Prop.~19.2,
p.~122]{Humphreys-alg}). Because maximal compact subgroups are conjugate
\see{cpct-conj}, there is some $v \in V$ with $T \subset N^v$. Because
$v$ normalizes~$C$ (indeed, it centralizes~$C$), and \ref{Ku-cent}
implies that $(X \times L)^\ka \subset C$, we know that $K^v$ contains
$(X \times L)^\ka$, so there is no harm in replacing~$K$ with~$K^v$.
Thus, we may assume $T \subset K$. Then, for any $r \in Q$, $(x,w)^\sa$
and $(x,w)^\ka$ are the projections of $(x,w)^\alpha$ into $W$~and~$T$,
respectively. Because $T$ centralizes~$W$, we conclude that $(X \times
Q)^\ka$ centralizes $W^\sa$, as desired.
 \end{proof}

\begin{lem} \label{wt-homo}
 Let $Q$ be a one-parameter subgroup of~$R$ that is normalized, but not
centralized, by~$A$. If $\sa|_Q$ is a crossed homomorphism, then
$\sa|_Q$ is a homomorphism.
 \end{lem}

 \begin{proof}
 For convenience, let $V$ be the Zariski closure of $\langle Q^\sa
\rangle$, $N = N_K(V)$ and $C = C_K(V)$. From \ref{crossed-K-homo}, we
know that the induced cocycle $\kq \colon X \times Q\to N/C$ is a
homomorphism. We wish to show that $\kq$ is trivial.

Because $A$ does not centralize~$Q$, there is some $a \in A$ with $r^a =
r^2$ for all $r \in Q$. Because $r^{a\sa} = r^{\sa (x,a)^\ka a^\ma}$, we
see that $(x,a)^\ka \in N$ \see{must-norm}, and that $(\cdot, a)^\ka$ is
constant, modulo~$C$. Thus, $a^\kq = (\cdot, a)^\ka C$ is a well-defined
element of $N/C$, so $\kq$ extends to a homomorphism defined on $QA$.
Thus, we have
 $r^{\kq a^\kq} = r^{a\kq} = (r^2)^\kq = (r^\kq)^2$, for all $r \in Q$.
 If $Q^\kq$ is nontrivial, this implies that $2$~is an eigenvalue of
$\Ad_{N/C} a^\kq$. But eigenvalues in a compact group all have absolute
value~$1$---contradiction.
 \end{proof}

\begin{lem}[{(cf.~\ref{Ku-cent})}] \label{K-cent} Given $r , s \in R$,
let $w = [s,r]\in R$. If
 $( \cdot ,r)^\sa = r^\sa$,
 $( \cdot ,s)^\sa = s^\sa$, and
 $( \cdot ,w)^\sa = w^\sa$ are essentially constant, and
$(X \times \langle s,w \rangle )^\ka$
centralizes~$ \langle s,w \rangle ^\sa$, then $(X \times s)^\ka$
centralizes~$r^\sa$.
 \end{lem}

\begin{proof}
 Same as \ref{Ku-cent} (and~\ref{u-equivariant}), with $s$~and~$s^\sa$
in place of $u$~and~$u^\ma$.
 \end{proof}

\begin{cor} \label{gen-crossed}
 Let $P$ and~$Q$ be subgroups of~$R$, such that $\sa|_P$ is a
homomorphism, $\sa|_Q$ is a crossed homomorphism, and $[P,Q] \subset
P$. Then $\sa|_{PQ}$ is a crossed homomorphism.
 \end{cor}

\begin{proof}
 Because $[P,Q] \subset P$, we see from \ref{K-cent} that $(X \times
P)^\ka$ centralizes~$Q^\sa$. Thus, for any $p \in P$ and $q \in Q$, we
have
 $$ ( \cdot ,pq)^\sa = p^\sa q^{\sa ( \cdot ,p)^{-\ka}} = p^\sa q^\sa $$
 is essentially constant, as desired.
 \end{proof}

\begin{cor} \label{gen-homo}
 Let $P$ and~$Q$ be subgroups of~$R$, such that $\sa|_P$ and $\sa|_Q$
are homomorphisms, and $[P,Q] \subset P$. Then $\sa|_{PQ}$ is a
homomorphism.
 \end{cor}

\begin{proof}
 From the preceding corollary, we know that $\sa|_{PQ}$ is a crossed
homomorphism.

Because $\sa|_P$ is a homomorphism, we know that $(X \times P)^\ka$
centralizes $P^\sa$. So, from \ref{K-cent} (and the fact that $[Q,P]
\subset P$), we see that $(X \times P)^\ka$ centralizes $Q^\sa$.

Because $P$ is nilpotent \see{special-case}, we may assume, by induction
on $\dim P$, that $\sa|_{Q[Q,P]}$ is a homomorphism, so $(X
\times Q)^\ka$ centralizes $[Q,P]^\sa$. Therefore, from \ref{K-cent} ,
we see that $(X \times Q)^\ka$ centralizes $(X \times P)^\sa$.

By combining the conclusions of the preceding two paragraphs, we
conclude that $(X \times QP)^\ka$ centralizes $(X \times QP)^\sa$, as
desired.
 \end{proof}

\begin{notation}
 Fix a normal subgroup~$Q$ of~$G$, contained in~$R$, such that
$\sa|_{[Q,Q]}$ is a homomorphism.
 \end{notation}

\begin{prop} \label{CQU-homo}
 We may assume $\sa|_{C_Q(U)}$ is a homomorphism.
 \end{prop}

\begin{proof}
 Let $\phi \colon X \to K_u$ be as in \ref{wolg-crossed}. From
\ref{CRL-homo} and~\ref{gen-homo}, we know that $\sa|_{C_R(L)[Q,Q]}$ is a
homomorphism. Thus, from~\ref{Kunip-cent-W}, we know that $K_u$
centralizes $C_R(L)[Q,Q]$, so $\sa^\phi|_{C_R(L)[Q,Q]} =
\sa|_{C_R(L)[Q,Q]}$ is a homomorphism. Thus, there is no harm in
replacing~$\alpha$ with~$\alpha ^\phi$, in which case, from the choice
of~$\phi$, we see that $\sa|_{C_Q(U)}$ is a crossed homomorphism. In
addition, \ref{wt-homo} (plus the fact that $\sa|_{C_Q(L)}$ is a
homomorphism) implies that $C_Q(U)$ is generated by subgroups~$T$ (one
of which is $[Q,Q]$), such that $\sa|_T$ is a homomorphism. Therefore,
\ref{gen-homo} implies that $\sa|_Q$ is a homomorphism.
 \end{proof}

For each $L$-module~$V$, we now define an
$AU$-submodule~$V^+$ and and $AU^-$-submodule~$V^-$. The specific
definition does not matter; what we need are the properties described in
the proposition that follows.

\begin{defn}
 Let $\Phi$ be the system of $\real$-roots of~$L$, let $\Delta$ be the
base for~$\Phi$ determined by~$U$, and let $\langle \Delta \rangle$ be
the $\integer$-span of~$\Delta$. (Note that $\Phi \subset \langle \Delta
\rangle$.) Define
 $$ \langle \Delta \rangle ^+
 = \set\{ \sum_{\alpha \in \Delta} k_\alpha \alpha \in \langle \Delta
\rangle \mid \exists\alpha, k_\alpha > 0 \} ,$$
 and $\langle \Delta \rangle ^-
 = - \langle \Delta \rangle ^+$.
 \end{defn}

\begin{defn}
 Let $V$ be a finite-dimensional $L$-module.
 For each linear functional on~$A$, we have the corresponding weight
space~$V_\lambda$. In particular, $V_0 = C_V(A)$.
 \begin{itemize}
 \item If $V$ is trivial, let $V^+ = V^- = V$.
 \item If $V$ is irreducible, and $V_0 = 0$, let $V^+ = V^- = V$.
 \item If $V$ is nontrivial and irreducible, and $V_0 \not= 0$, then
every weight of~$V$ belongs to $\langle\Delta \rangle$; let
 $V^+ = \sum_{\lambda \in \langle \Phi \rangle ^+} V_\lambda$ and
 $V^- = \sum_{\lambda \in \langle \Phi \rangle ^-} V_\lambda$.
 \item In general, define $V^+ = \sum W^+$ and $V^- = \sum W^-$, where
each sum is over all irreducible submodules~$W$ of~$V$.
 \end{itemize}
 \end{defn}

\begin{prop} \label{vpm-prop}
 Let $V$ be a finite-dimensional $L$-module. Then:
 \begin{enumerate}
 \item \label{vpm-span}
 $V = V^- + V_0 + V^+$;
 \item \label{vpm-inter}
 $V^+ \cap V_0 = V^- \cap V_0 = C_V(L)$;
 \item \label{vpm-submodule}
 $V^+$ is an $AU$-submodule of~$V$;
 \item
 $V^-$ is an $AU^-$ submodule; and
 \item \label{vpm-gen} $V_0 + V^+$ is the smallest $U$-submodule of~$V$
that contains both $V_0$ and $V^+ \cap V^-$.
 \end{enumerate}
 \end{prop}

\begin{proof}
 Only \pref{vpm-gen} is perhaps not clear from the definition. We may
assume $V$ is irreducible. Assume, furthermore, that $V$ is nontrivial
and $V_0 \not= 0$, for otherwise we have $V^+ \cap V^- = V$, so the
desired conclusion is obvious. Fix some $\lambda \in \langle \Delta
\rangle ^+$. If $\lambda \not\in \langle \Delta \rangle^-$, then, in the
unique representation $\lambda = \sum k_\alpha \alpha$ of~$\lambda$ as a
linear combination of the elements of~$\Delta$, it must be the case that
every~$k_\alpha$ is nonnegative. Because $||\lambda || > 0$, this
implies that there is some $\beta \in\Delta$ whose inner product
with~$\lambda$ is strictly positive. Let $\Lie U$, $\Lie A$, and~$\Lie
L$ be the Lie algebras of $U$, $A$, and~$L$, respectively. By induction
on $\sum k_\alpha$, we may assume that $V_{\lambda -\beta }$ is in the
$\Lie U$-submodule of~$V$ generated by $V_0$ and $V^+ \cap V^-$, so it
suffices to show that $[\Lie L_\beta ,V_{\lambda -\beta }] = V_\lambda$. 

Choose $u \in \Lie L_\beta$, $h \in \Lie A$, and $v \in \Lie L_{-\beta
}$, such that the linear span $\langle u,h,v \rangle$ is a subalgebra
of~$\Lie L$ isomorphic to $\sl_2(\real)$, and such that $h^\mu$~is equal
to the inner product of~$\beta$ with~$\mu$, for every weight~$\mu$ (see
\cite[Eqn.~(7) of \S IX.1, p.~407]{Helgason}). Then, if we restrict~$V$
to a representation of $\langle u,h,v \rangle$, we see that vectors in
the space~$V_\lambda$ have strictly positive weight, so the structure
theory of $\sl_2(\real)$-modules implies that $[u, V_{\lambda -\beta }]
= V_\lambda $.
 \end{proof}

\begin{defn}
 Let $\Lie Q$ be the Lie algebra of~$Q$. Clearly, $\Lie Q^+ [\Lie Q,
\Lie Q]$ and $\Lie Q^+ [\Lie Q, \Lie Q]$ are Lie subalgebras of~$\Lie
Q$; let $Q^+$ and~$Q^-$ be the corresponding connected Lie subgroups
of~$Q$.

In addition, we let $Q_0 = C_Q(A) [Q,Q]$ and $Q_0^+ = Q_0 Q^+$.

We define $Q^-$ and $Q_0^-$ analogously.

Let $Q^\pm = Q^+ \cap Q^-$, and $Q_0^\pm = Q_0 Q^\pm$.
 \end{defn}

\begin{prop}
 We may assume $\sa|_{Q^+}$ is a homomorphism, and $\sa|_{Q_0^+}$ is a crossed
homomorphism.
 \end{prop}

\begin{proof}
 ($Q^+$) For a proof by induction, it suffices to show that if $T$ is a
one-parameter subgroup of~$Q$ that is normalized by~$A$, and $N$ is a
subgroup of~$Q^+$, normalized by~$TU$, such that $[T,U] \subset N$ and
$\sa|_N$ is a homomorphism, then we may assume $\sa|_{TN}$ is a
homomorphism. From \ref{wolg-const}, we see that there is a Borel
function $\phi \subset X \to K_u$, such that $\sa^\phi |_T$ is a crossed
homomorphism. From~\ref{Kunip-cent-W}, we see that $K_u$
centralizes~$N^\sa$, so there is no harm in replacing~$\alpha$
with~$\alpha ^\phi$, so we may assume $\sa|_T$ is a crossed
homomorphism.

If $T \subset C_Q(A)$, then, by definition of~$Q^+$
(and~\ref{vpm-prop}\pref{vpm-inter}), we must have $T \subset
C_Q(L)[Q,Q]$, so $\sa|_T$ is a homomorphism (see~\ref{CRL-homo}
and~\ref{gen-homo}). On the other hand, if $T \not\subset C_Q(A)$, then
\ref{wt-homo} applies, so we again conclude that $\sa|_T$ is a
homomorphism. Thus, from \ref{gen-homo}, we conclude that $\sa|_{TN}$ is
a homomorphism, as desired.

($Q_0^+$) From \ref{wolg-const} (and~\ref{Kunip-cent-W}), we see that we
may assume $\sa|_{Q_0}$ is a crossed homomorphism. So \ref{gen-crossed}
implies that $\sa_{Q_0^+}$ is a crossed homomorphism.
 \end{proof}

Similarly, by considering the opposite unipotent subgroup~$U^-$, we
obtain:

\begin{prop} \label{Q-minus}
 There is a Borel function $\phi \colon X \to K$, such that the
restriction $\sa^\phi |_{Q^- C_R(L)}$ is a homomorphism, and $\sa^\phi
|_{Q_0^- C_R(L)}$ is a crossed homomorphism. \qed
 \end{prop}

\begin{prop} \label{wma-crossed}
 We may assume $\sa|_{Q^-}$ is a homomorphism, and $\sa|_Q$ is a crossed
homomorphism.
 \end{prop}

 \begin{proof} Because $Q = Q^- Q_0^+$ and $[Q,Q] \subset Q^-$, we see
from \ref{gen-crossed} that the second conclusion follows from the
first. Choose a Borel function $\phi \colon X \to K$, as in
Prop.~\ref{Q-minus}. It suffices to show that we may assume $X^\phi =e$.

\step{1}{We may assume $X^\phi$ centralizes $(Q_0^\pm)^\sa$.}
 For any $r \in Q_0^\pm$ and a.e.~$x \in X$, we have $r^{\sa^\phi } =
r^{\sa x^\phi }$, so there is some $k \in K$, such that, for a.e.~$x \in
X$, we have
 $x^\phi \in C_K\bigl((Q_0^\pm)^\sa\bigr) k$. We may replace $x^\phi$
with the function $x \mapsto x^\phi k^{-1}$.

\step{2}{$X^\phi$ centralizes $(Q_0^\pm)^{\sa k}$, for every $k \in
(X \times Q_0^\pm)^\ka$.}
 For
$r,s \in Q_0^\pm$, and a.e.~$x \in X$, we have
 $(rs)^\sa = r^\sa s^{\sa (x,r)^{-\ka}}$. Because
$X^\phi$ centralizes $(rs)^\sa$ and~$r^\sa$ (see Step~1), this implies
that $X^\phi$ centralizes $s^{\sa (x,r)^{-\ka}}$, as desired.

\step{3}{For every $r \in Q_0^+$ and
 $k \in \langle (X \times Q_0^+)^\ka \rangle$, there is some
 $k' \in (X \times Q_0)^\ka$ with $r^k = r^{k'}$.}
 First note that, modulo $C_K\bigl((Q_0^+)^\sa\bigr)$, the
cocycle~$\ka|_{Q_0^+}$ is a homomorphism \see{crossed-K-homo}. Because
the image of a homomorphism is always a subgroup, this implies that, for
any $k \in \langle (X \times Q_0^+)^\ka \rangle$, there is some $r \in
Q_0^+$ with
 $( \cdot ,r)^\ka \in C_K\bigl((Q_0^+)^\sa\bigr) k$~a.e.

Therefore, it will suffice to show that $(X \times Q^+)^\ka$
centralizes~$(Q_0^+)^\sa$. Because $\sa|_{Q^+}$ is a homomorphism (so $(X
\times Q^+)^\ka$ centralizes $(Q^+)^\sa)$, and $[Q_0^+, Q_0^+] \subset
Q^+$, Lem.~\ref{K-cent} provides this conclusion.

\step{4}{$(Q_0^+)^\sa$ is centralized by~$X^\phi $.}
 For $r \in Q_0^\pm$ and $u \in U$, we have
 $$
 (xu, r^u)^\alpha
 = (x,u)^{-\alpha } (x,r)^\alpha (xr,u)^\alpha
 \in (x,r)^{\alpha (x,u)^\ka u^\ma} K ,$$
 so
 $r^{u\sa} = r^{\sa (x,u)^\ka u^\ma} = r^{\sa u^\ma}$,
 because $(x,u)^\ka$ centralizes~$(Q_0^\pm)^\sa$ \see{Ku-cent}.
 Therefore, for $k \in \langle (X \times Q_0^\pm)^\ka \rangle$, we have
 $r^{u \sa k} = r^{\sa u^\ma k} = r^{\sa k u^\ma}$.
 Since $X^\phi$ centralizes~$u^\ma$ and $r^{\sa k}$ (for the
latter, see Steps~2 and~3), this implies that $X^\phi$ centralizes~$r^{u
\sa k}$. Thus, because
 \begin{equation*}
 \begin{split}
 (Q_0^+)^\sa
 &= \langle r^u \mid r \in Q_0^\pm, u
\in U \rangle ^\sa 
 \qquad\qquad \textnormal{(see \ref{vpm-prop}\pref{vpm-gen})} \\
 &\subset
 \bigl\langle r^{u\sa k}
 \mid r \in Q_0^\pm, u \in U, k \in \langle (X \times Q_0^\pm)^\ka
\rangle \bigr\rangle , \\
 \end{split}
 \end{equation*}
 we have the desired conclusion.

\step{5}{We may assume $X^\phi =e$.}
 From Step~4, we see that replacing $\alpha$ with the equivalent cocycle
$\alpha ^\phi$ will not change $\sa|_{Q_0^+}$, so we may assume $x^\phi
=e$.
 \end{proof}

\begin{prop} \label{Q-homo}
 $\sa|_Q$ is a homomorphism.
 \end{prop}

\begin{proof}
 For $g \in G$, let $(x,g)^\pa = (x,g)^\sa (x,g)^\ma$.
 For $r \in Q$ and $u \in U$, because $( \cdot ,r)^\ka$
centralizes~$u^\ma$, we have
 $( \cdot , r u) ^\pa = r^\sa u^\ma$
 is essentially constant. Thus, $\pa|_{QU}$ is a crossed homomorphism,
so the cocycle $\ka|_{QU}$ is a homomorphism, modulo
$C_K\bigl((QU)^\sa\bigr)$ \see{crossed-K-homo}. Then, because $QU$ is
solvable, the fact that connected, compact, solvable groups are abelian
\see{cpct-abel}, implies that $[QU,QU]^\ka \subset
C_K\bigl((QU)^\sa\bigr)$, so $\pa|_{[QU,QU]}$ is a homomorphism. In
particular, $\sa|_{Q_0 \cap [U,Q]}$ is a homomorphism.

Now $Q_0$ is generated by $Q_0 \cap [U,Q]$, $C_Q(L)$, and $[Q,Q]$. The
restriction of~$\sa$ to each of these subgroups is a homomorphism, so we
conclude from~\ref{gen-homo} that $\sa|_{Q_0}$ is a homomorphism. Then,
because $Q$ is generated by $Q_0$, $Q^-$, and~$Q^+$, we conclude
from~\ref{gen-homo} that $\sa|_Q$ is a homomorphism.
 \end{proof}

This proposition provides the induction step in a proof of the following
important corollary:

\begin{cor}
 We may assume $\sa|_R$ is a homomorphism. \qed
 \end{cor}

\begin{prop} \label{Kcent-triv}
 We have $C_K( R^\sa ) = e$.
 \end{prop}

\begin{proof} Because $C_K(S)$ is centralized (hence normalized) by~$M$
and~$S$, and is normalized by~$K$, we see that $C_K(S)$ is a normal
subgroup of~$H$. Being also compact, it must be trivial. Thus, letting
$S'$ be the Zariski closure of $\langle R^\sa \rangle$, it suffices to
show $S' \normal H$, for then $S' = S$. Furthermore, because $G = LR$, we
need only show that $(X \times L)^\alpha$ and $(X \times R)^\alpha$
each normalize~$S'$.

For $u \in L$ and $r \in R$, we have $r^{\sa (x,u)^\alpha} = r^{u\sa}
\in R^\sa$, so we see that $(x,u)^\alpha$ normalizes~$S'$.

 Because $\sa|_R$ is a homomorphism, we know that $(X \times R)^\ka$
centralizes~$S'$. It is obvious from the definition of~$S'$ that
$R^\sa$ normalizes~$S'$.
 \end{proof}

\begin{cor} \label{alpha-homo}
 $\alpha$ is a homomorphism.
 \end{cor}

\begin{proof}
 Because $\ma$ and~$\sa|_R$ are homomorphisms, it only remains to show
that $\ka$~is a homomorphism.

 For $u \in L$ and $r \in R$, we have $r^{u\sa} = r^{\sa (x,u)^\alpha}$,
so the proposition implies that $( \cdot ,u)^\alpha$ is essentially
constant.

Because $\sa|_R$ is a homomorphism, we know that $(X \times R)^\ka$
centralizes~$R^\sa$. Hence, the proposition implies that $(X \times
R)^\ka = e$.
 \end{proof}

\section{The remaining proofs} \label{remaining}

In this section, we describe how to finish the proof of
Thm.~\ref{mainthm} \see{general-case} and we prove Cor.~\ref{main-cor}
\see{main-cor-pf}.

\subsection{Proof of the remaining case of Theorem~\ref{mainthm}}
\label{general-case}

The notation and hypotheses of Thm.~\ref{mainthm} are in effect
throughout this subsection, but, unlike in~\S\ref{main-proof}, we do not
assume that $\Rad G$ is nilpotent, nor that $G$ has no nontrivial,
compact, semisimple quotients \see{special-case}. Let
 \begin{itemize}
 \item $N = \nil G$; and
 \item $R$ be the product of $\Rad G$ with the maximal compact factor of
a Levi subgroup of~$G$.
 \end{itemize}

\begin{prop}
 We may assume $\sa|_R$ is a crossed homomorphism.
 \end{prop}

\begin{proof}
 From the main proof (\S\ref{main-proof}), applied to the group $LN$, we
see that we may assume $\sa|_N$ is a homomorphism. Let $K_u$ be the
closure of $\langle (X \times U)^\ka \rangle$. Then \ref{wolg-crossed}
implies there is a Borel function $\phi \colon X \to K_u$, such that
$\sa^\phi|_{C_R(U)}$ is a crossed homomorphism. Because $K_u$
centralizes~$N^\sa$ \see{Ku-cent}, we have $\sa^\phi |_N = \sa|_N$, so
there is no harm in replacing $\alpha$ with the cohomologous
cocycle~$\alpha ^\phi $. Thus, we may assume $\sa|_{C_R(U)}$ is a
crossed homomorphism. Then, because $R = N \, C_R(U)$ and $[N,R] \subset
N$, Lem.~\ref{gen-crossed} implies that $\sa|_R$ is a crossed
homomorphism.
 \end{proof}

\begin{prop}[\cf{Kcent-triv}]
 We have $C_K( R^\sa ) = e$.
 \end{prop}

\begin{proof}
 Let $S'$ be the Zariski closure of $\langle R^\sa \rangle$. As in the
proof of~\ref{Kcent-triv}, we wish to show $S'$ is normalized by $(X
\times L)^\alpha$ and $(X \times R)^\alpha$.

For $u \in L$ and $r \in R$, we have $r^{u\sa} = r^{\sa (x,u)^\alpha}$,
so we see that $(x,u)^\alpha$ normalizes~$S'$.

For any $r,s \in R$, we have
 $(rs)^\sa = r^\sa s^{\sa (x,r)^{-\ka}}$, so $s^{\sa (x,r)^{-\ka}} \in
\langle R^\sa \rangle$. This implies that $(X \times R)^\ka$
normalizes~$S'$. It is obvious from the definition of~$S'$ that $R^\sa$
also normalizes~$S'$.
 \end{proof}

\begin{cor}[\cf{alpha-homo}]
 $\alpha$ is a homomorphism.
 \end{cor}

\begin{proof}
 Because $G = LR$, it suffices to show that $\alpha |_L$ and $\alpha
|_R$ are homomorphisms.

 For $u \in L$ and $r \in R$, we have $r^{u\sa} = r^{\sa (x,u)^\alpha}$,
so the proposition implies that $( \cdot ,u)^\alpha$ is essentially
constant. Therefore, $\alpha |_L$ is a homomorphism.

For any $r,s \in R$, we have
 $(rs)^\sa = r^\sa s^{\sa (x,r)^{-\ka}}$, so the proposition implies that
$( \cdot ,r)^\ka$ is essentially constant. Therefore, $\ka |_R$ is a
homomorphism. Because $\sa|_R$ is a crossed homomorphism, this implies
that $\alpha |_R$ is a homomorphism.
 \end{proof}

\subsection{Proof of Corollary~\ref{main-cor}} \label{main-cor-pf}

The notation and hypotheses of Cor.~\ref{main-cor} are in effect
throughout this subsection. We wish to verify the hypotheses of
Thm.~\ref{mainthm}.

By passing to an ergodic component of $X \mathbin{ \times _\alpha }
H/H^\circ$, which is a finite cover of~$X$, we may assume $H$ is
connected (cf.~\cite[Prop.~9.2.6, p.~168]{ZimmerBook}).
Then we
may write $H = S \prodsemi (MK)$, where
 \begin{itemize}
 \item $S$ is a connected, split, solvable subgroup;
 \item $M$ is connected and semisimple, with no compact
factors;
 \item $K$ is a compact subgroup that centralizes~$M$; and
 \item $M \cap K = Z(M)$.
 \end{itemize}

By assumption, the Zariski hull of~$\alpha|_{X \times L}$ is reductive.
Then, because $L$ has the Kazhdan Property (see \cite[Thm.~7.1.4,
p.~130]{ZimmerBook}), we see that the center of this Zariski hull is
compact (see \cite[Thm.~9.1.1, p.~162]{ZimmerBook}). Thus, the Zariski
hull is contained in (a conjugate of) $MK$, so, by replacing~$\alpha $
with an equivalent cocycle, we may assume $(X \times L)^\alpha \subset
MK$.

Now $\alpha$ induces a cocycle $\aq \colon X \times G \to H/(SK) \iso
M/(M \cap K)$. 

Although the statement of Zimmer's Theorem \pref{Zimmer-super} assumes
that $G$ is semisimple, the proof shows that this hypothesis is not
necessary if the Zariski hull of the cocycle is assumed to be
semisimple. Thus, we have the following:

\begin{thm}[{(Zimmer, cf.~pf.~of \cite[Thm.~5.2.5, p.~98]{ZimmerBook})}]
 $\aq$ is cohomologous to a homomorphism. \qed
 \end{thm}

So there is no harm in assuming that $\bar\alpha$ itself is a
homomorphism. By replacing $G$ with a finite cover \see{cover-coho}, we
may assume that $\aq$ lifts to a homomorphism $\sigma \colon G \to M$. 

Note that, because $(X \times L)^\alpha \subset MK$, the definition
of~$\sigma$ implies that, for every $l \in L$, we have $(x,l)^\alpha
\in l^\sigma K$ for a.e.~$x \in X$. 

Let $J$ be the product of $\Rad G$ with the maximal compact factor of a
Levi subgroup of~$G$; then $J$ is a connected, normal subgroup of~$G$.
Therefore, because $G^\aq$ is Zariski dense in $M/(M \cap K)$, we see
that $J^\aq$ is normal in~$M/(M \cap K)$. But $M/(M \cap K)$ has no
nontrivial, normal subgroups that are solvable or compact, so
we conclude that $J^\aq$ is trivial. Because $G = LJ$, this implies that
$L^\sigma$ is Zariski dense in~$M$. Because connected, semisimple
subgroups have finite index in their Zariski closure
\cite[Thm.~VIII.3.2, p.~112]{Hochschild-alg}, this implies $L^\sigma =M$,
so $H = S L^\sigma K$.

Thus, the hypotheses of Thm.~\ref{mainthm} are verified, so we conclude
that $\alpha$ is cohomologous to a homomorphism. This completes the
proof of Cor.~\ref{main-cor}. \qed

\end{document}